\documentclass{amsart}

\usepackage{graphicx}
\usepackage{longtable}
\usepackage{amsmath}
\usepackage{amssymb}
\usepackage[pdfauthor={Richard J. Mathar},colorlinks=true,citecolor=blue]{hyperref}

\newtheorem{rem}{Remark}
\newtheorem{alg}{Algorithm}
\newtheorem{exa}{Example}

\def\deg{\mathop{\mathrm{deg}}\nolimits}

\numberwithin{equation}{section}

\begin{document}

\title[Series of Generalized Fresnel Integrals]{Series Expansion of Generalized Fresnel Integrals}

\author{Richard J. Mathar}
\urladdr{http://www.mpia.de/~mathar}
\email{mathar@mpia.de}
\address{Leiden Observatory, Leiden University, P.O. Box 9513, 2300 RA Leiden, The Netherlands}


\date{\today}

\begin{abstract}
The two Fresnel Integrals are real and imaginary part of the integral over complex-valued
$\exp(ix^2)$ as a function of the upper limit. They are special cases of the integrals over $x^m\exp(ix^n)$
for integer powers $m$ and $n$, which are essentially Incomplete Gamma Functions. We generalize
one step further
and focus on evaluation of the integrals with kernel $p(x)\exp[i\phi(x)]$ and polynomials $p$ and $\phi$.
Series reversion of $\phi$ seems not to
help much,
but repeated partial integration 
leads to a first order differential equation for an auxiliary oscillating function which allows
to fuse the integrals and their complementary integrals.
\end{abstract}

\maketitle

\section{Introduction}
\subsection{Motivation}

A criticism to the orthogonal basis of the Zernike Polynomials on the unit circle
is the lack of the minimax property--if compared for example to the Chebyshev
Polynomials on the one-dimensional interval \cite{FraserJACM12}. The radial polynomials of higher order,
starting already with $Z_3^{\pm 1}$, $Z_4^0$ and $Z_4^{\pm 2}$, attain amplitudes
at the local maxima and minima within the circle which are smaller (in absolute
value) than the amplitude at the rim.
One way to equilibrate these, to enhance the extrema in the interior
and to weaken the extrema along the perimeter, is to add higher powers of the radial distance
to the polynomials. An alternative, which triggered this work, may employ sines and cosines
of polynomials of the radial distance which evidently oscillate up and down
with constant amplitude. The requirement of orthogonality and normalization then
puts focus on radial integrals of the form
\begin{equation}
I_{p,\phi}(u)=\int_0^u p(x)e^{i\phi(x)}dx,
\label{eq.Idef}
\end{equation}
where $\phi$ is the sum
or difference of two phase-polynomials of the basis.

The concept applies to Zernike Polynomials in higher dimensions as well, so
$p$ are (small) integer powers of the radial coordinate, representing the
radial part of the Jacobian determinant of the transformation from Cartesian to Hyper-spherical Coordinates. 
(We assume that the angular variables are represented by Hyperspherical Harmonics,
which ensures that these oscillations oscillate with amplitudes of equal magnitude
already by design.)
This work deals with the evaluation of the integrals $I_{p,\phi}$ given 
a finite upper limit $u$ plus the polynomials $p$ and $\phi$. The
adjustment (variation) of the coefficients of a set of polynomials $\phi$ to achieve orthogonality
and  to establish an explicit basis of this type is not considered here.

\subsection{Overview}
Section \ref{sec.simpl} summarizes some elementary closed-form solutions.
Section \ref{sec.fres} is an overview of known results concerning
the cases with constant $p(x)$ and $\phi(x)$ being a power of $x$, which includes
the Fresnel Integrals \emph{eo ipso}.
Section \ref{sec.powp} widens the perspective to the cases where
both $p(x)$ and $\phi(x)$ are integer powers of $x$. The evaluation is covered by
the standard theory of the Incomplete Gamma Function. 
Section \ref{sec.rev} looks at series reversions if $\phi(x)$ are polynomials,
leading to barren results.
Section \ref{sec.pract} shows that asymptotic expression for the complementary
integral with one limit bound at infinity is obtained easily via an auxiliary
inhomogeneous differential equation.
The fusion of the series representations of that auxiliary function
at small upper limits on one hand
and its asymptotic series on the other is prosecuted by finding the
adaptive free parameter
of the general solution of the differential equation.

\subsection{Simple Cases}\label{sec.simpl}
The simplest cases are polynomials $\phi$ of zero or first degree \cite[2.635]{GR}:
\begin{eqnarray}
\int x \cos \alpha_0 dx &=& \frac{x^2}{2}\cos \alpha_0,\\
\int x \sin \alpha_0 dx &=& \frac{x^2}{2}\sin \alpha_0,\\
\int x \cos (\alpha_0 +\alpha_1 x)dx &=& \frac{1}{a_1^2}\cos (\alpha_0+\alpha_1 x)+\frac{x}{a_1}\sin(\alpha_0+\alpha_1x),\\
\int x \sin (\alpha_0+\alpha_1 x) dx &=& \frac{1}{a_1^2}\sin (\alpha_0+\alpha_1x)-\frac{x}{\alpha_1}\cos(\alpha_0+\alpha_1x).
\end{eqnarray}
Repeated partial integration yields
\begin{equation}
\int x^n e^{iax} dx
=
\frac{e^{iax}}{i a^{n+1}}\sum_{k=0}^n i^k \frac{n!}{(n-k)!} (ax)^{n-k}
,
\label{eq.xnexpax}
\end{equation}
with real and imaginary parts \cite[2.323,2.633]{GR}
\begin{equation}
\int x^n\cos(ax) dx=
\sum_{k=0}^n \frac{n!}{(n-k)!} \frac{x^{n-k}}{a^{k+1}}
\sin(ax+k\pi/2),
\end{equation}
\begin{equation}
\int x^n\sin(ax) dx=
-\sum_{k=0}^n \frac{n!}{(n-k)!} \frac{x^{n-k}}{a^{k+1}}
\cos(ax+k\pi/2).
\label{eq.xnsinax}
\end{equation}
(Erratum:
A factor $(-1)^k$ is missing in the last term of \cite[2.633.6]{GR}).
Further obvious special cases are
\begin{equation}
\int x^{n-1}\sin(x^n) dx = -\frac{1}{n}\cos(x^n);
\quad
\int x^{n-1}\cos(x^n) dx = \frac{1}{n}\sin(x^n) ;
\end{equation}
\begin{equation}
\int x^{n-1}e^{ix^n} dx = \frac{1}{in}e^{ix^n} .
\label{eq.triv}
\end{equation}

\begin{rem}\label{rem.o0}
The polynomial in the argument of the trigonometric function may
be assumed to have minimum low degree 1, because its constant term is easily
moved to the front of the integral with \cite[4.3.16,4.3.17]{AS}
\begin{equation}
\sin(\sum_{j=0}^l\alpha_j x^j) =
\sin \alpha_0 \cos (\sum_{j=1}^l\alpha_j x^j)
+\cos \alpha_0 \sin (\sum_{j=1}^l\alpha_j x^j),
\label{eq.sina0}
\end{equation}
\begin{equation}
\cos(\sum_{j=0}^l\alpha_j x^j) =
\cos \alpha_0 \cos (\sum_{j=1}^l\alpha_j x^j)
-\sin \alpha_0 \sin (\sum_{j=1}^l\alpha_j x^j),
\label{eq.cosa0}
\end{equation}
which are real and imaginary part of
\begin{equation}
\exp(i\sum_{j=0}^l\alpha_j x^j) =
\exp(i\alpha_0)
\exp(i\sum_{j=1}^l\alpha_j x^j)
.
\label{eq.phiofalpha}
\end{equation}
\end{rem}

Combinations of an odd polynomial $p(x)$ with an even polynomial $\phi(x)$
are reduced in a first step
by the substitution $y=x^2$,
\begin{equation}
\int_0^u x^{2m+1}\exp(ix^{2n}) dx = \frac{1}{2}\int_0^{u^2}  y^m\exp(iy^n) dy.
\label{eq.eodd}
\end{equation}

\section{Fresnel Integrals}\label{sec.fres}

\begin{figure}
\includegraphics[width=0.8\textwidth]{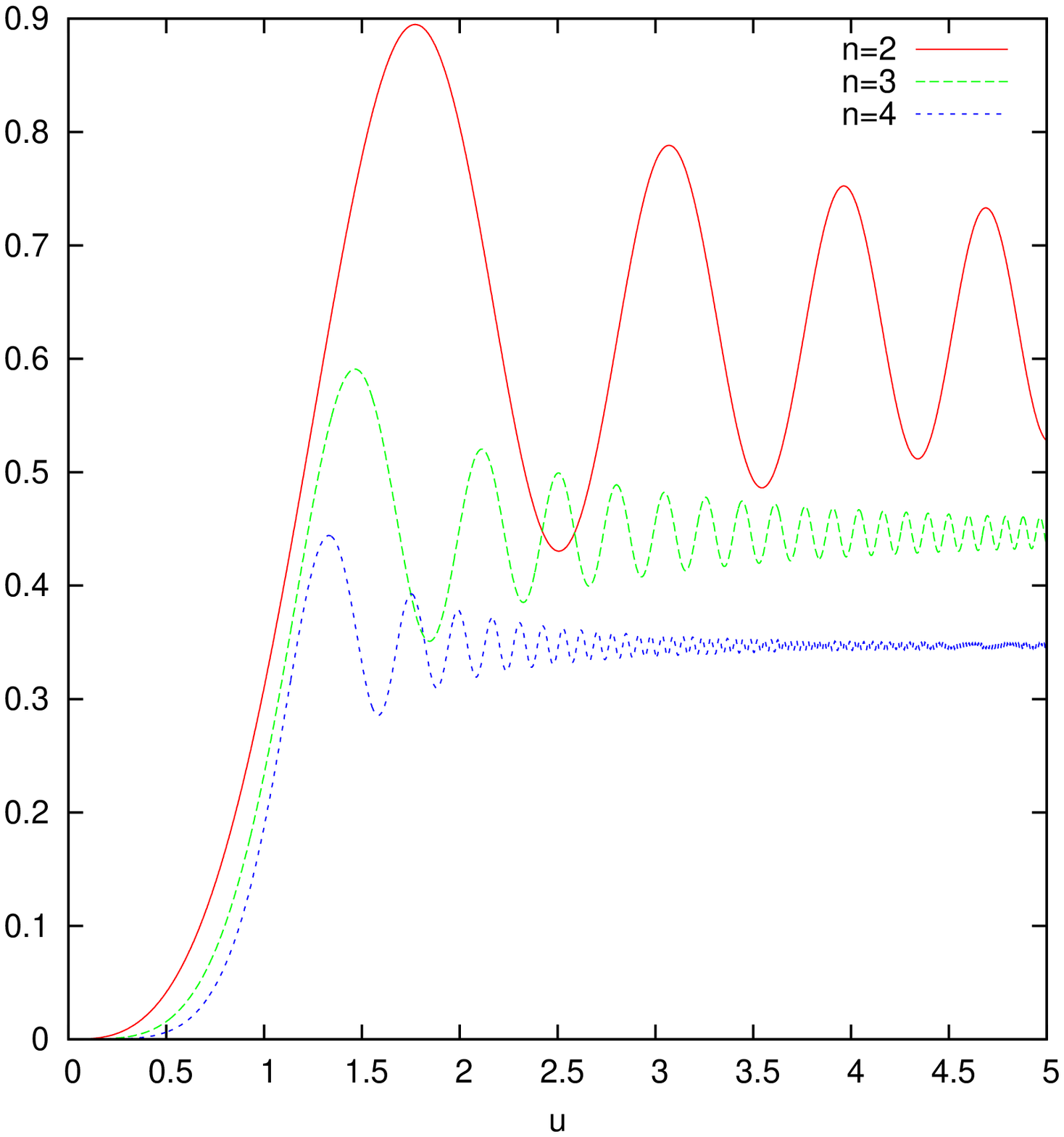}
\caption{The integral $\int_0^u \sin(x^n)dx$ for three different powers $n$.}
Extrema are at
$u^n=l\pi$ with $l=0,1,2,\ldots$.
\end{figure}
\begin{figure}
\includegraphics[width=0.8\textwidth]{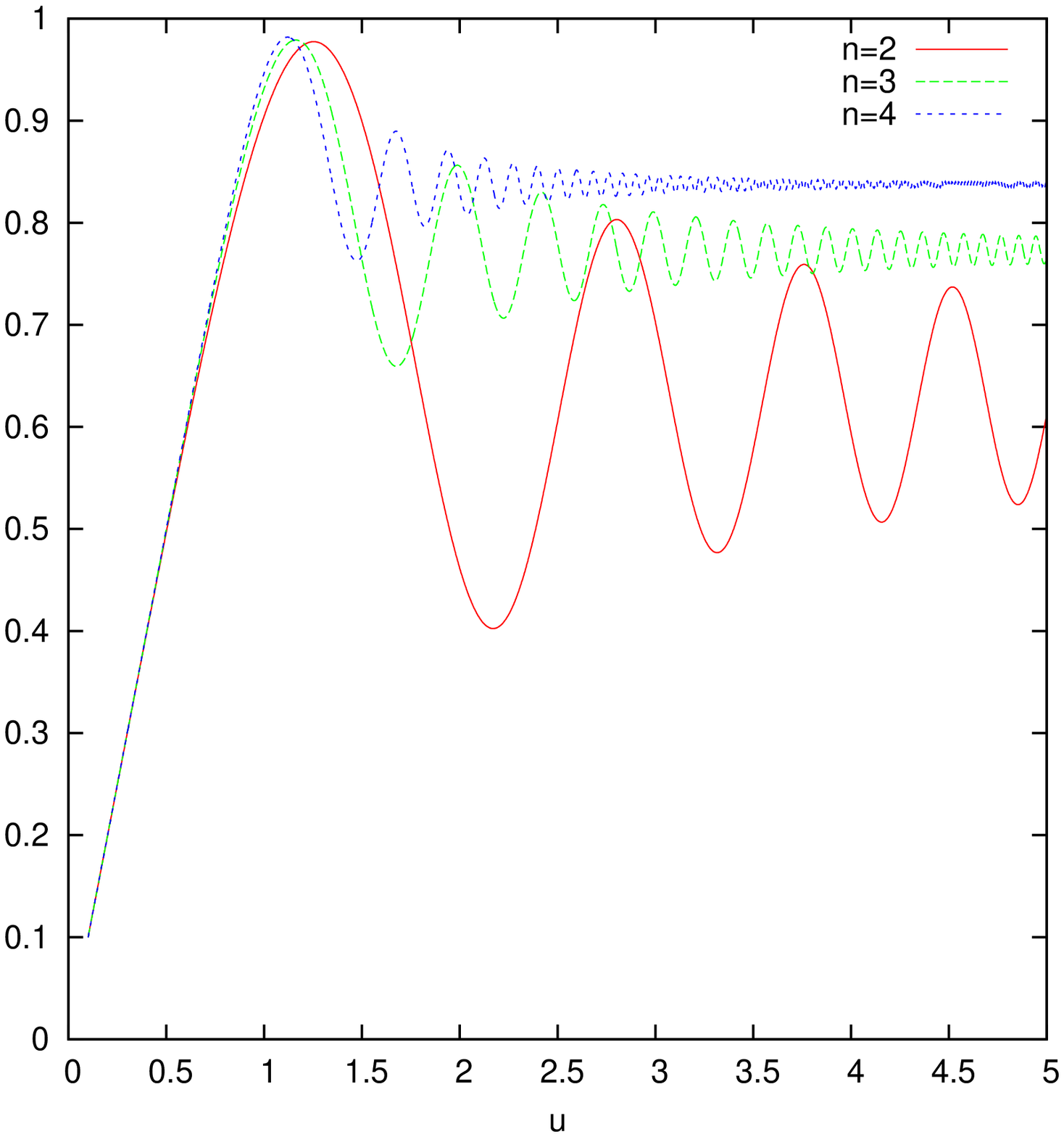}
\caption{The integral $\int_0^u \cos(x^n)dx$ for three different powers $n$.}
Extrema are at
$u^n=(l+1/2)\pi$ with $l=0,1,2,\ldots$.
\end{figure}
\begin{figure}
\includegraphics[width=0.8\textwidth]{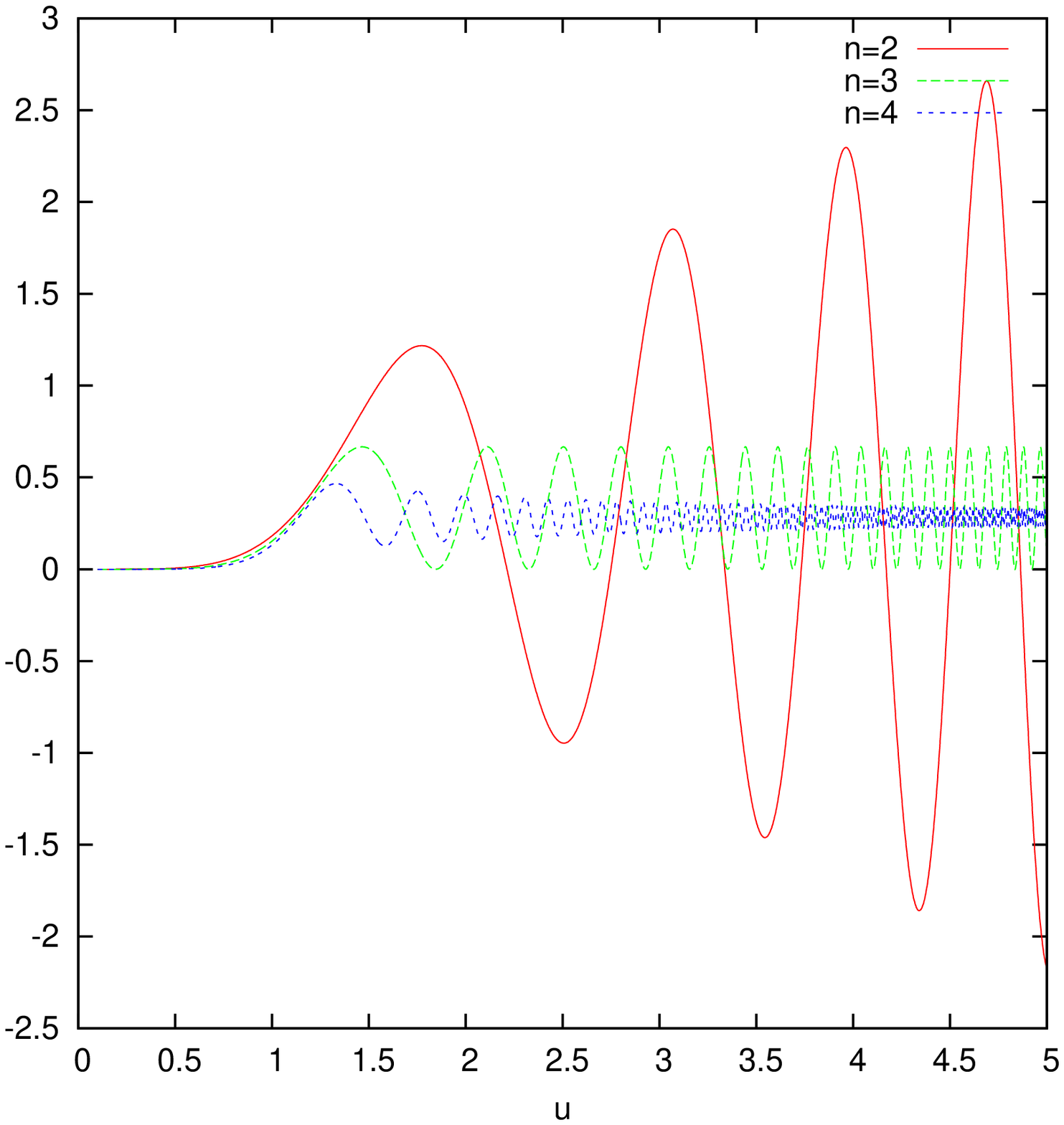}
\caption{The integral $\int_0^u x^2 \sin(x^n)dx$ for three different values of $n$.}
\label{fig.Jn22to4S}
\end{figure}
\begin{figure}
\includegraphics[width=0.8\textwidth]{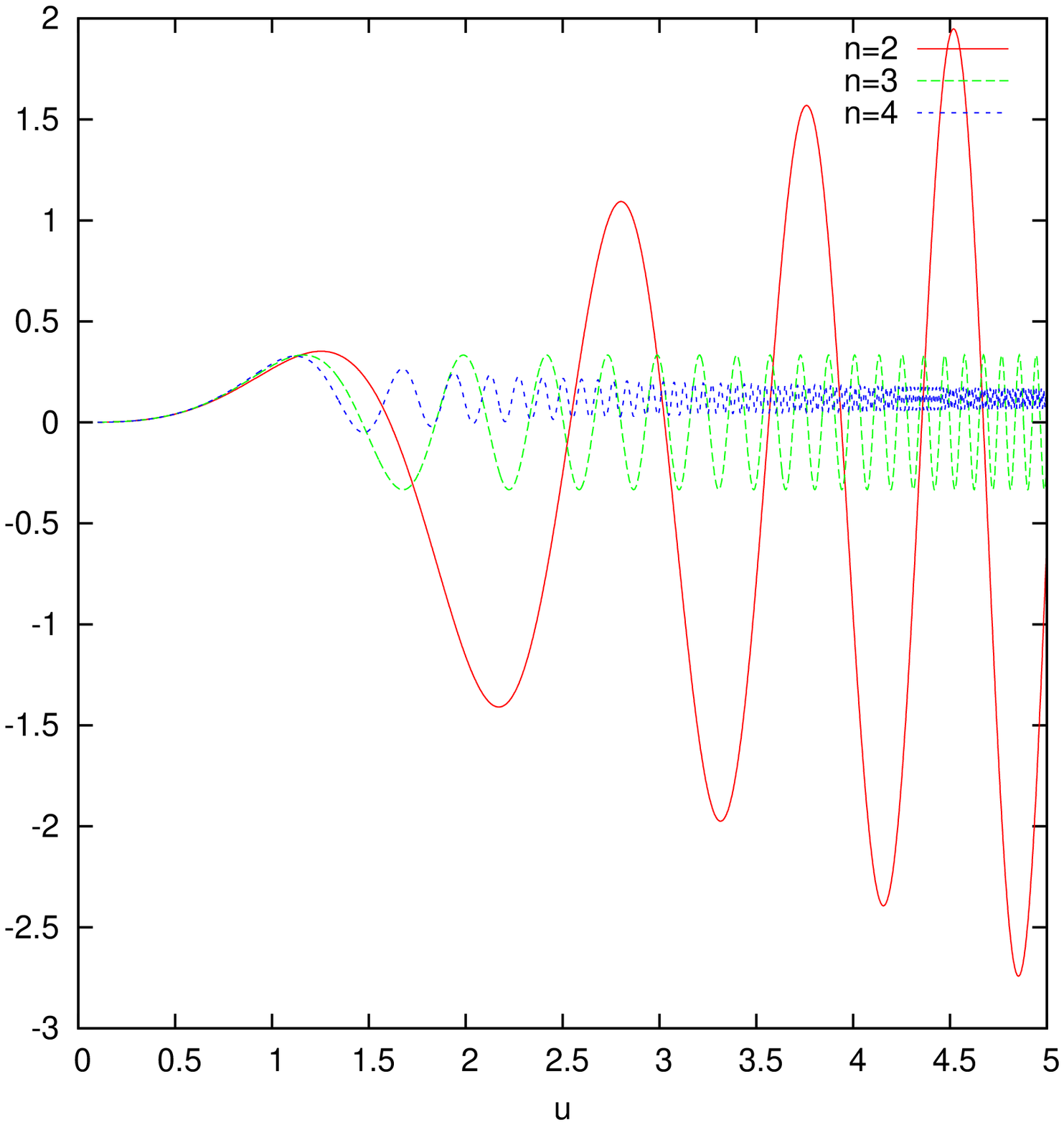}
\caption{The integral $\int_0^u x^2\cos(x^n)dx$ for three different values of $n$.}
\label{fig.Jn22to4C}
\end{figure}
The integrals $\int_0^u \sin(x^2)dx$ and $\int_0^u  \cos(x^2)dx$
define---up to some unimportant constant scale factor---the Fresnel
Integrals \cite[7.3]{AS}\cite[8.25]{GR}\cite{BoersmaMCom14,HangelJEM1,Wrenchdtic,BulirschNumMath9,vanWijngaardKNAW}.
In this section we look at the generalized integrals
$\int\sin(x^n)dx$ and $\int \cos(x^n)dx$ for integer $n=2,3,4,\ldots$.

The substitution $x=y^{1/n}$ followed by iterated partial integration generates
\begin{multline}
\int e^{ix^n}dx = \frac{1}{n}\int y^{1/n-1}e^{iy} dy
\\
=
\frac{1}{in} y^{1/n-1}e^{iy}
-\frac{1}{in}(\frac{1}{n}-1)\int y^{1/n-2}e^{iy} dy
\\
=
\frac{1}{in} y^{1/n-1}e^{iy}
-\frac{1}{in}(\frac{1}{n}-1)
-\frac{1}{i^2n}(\frac{1}{n}-1)
y^{1/n-2}e^{iy}
+\frac{1}{i^2n}(\frac{1}{n}-1)
(\frac{1}{n}-2)\int y^{1/n-3}e^{iy} dy
=\ldots
\end{multline}
The product of the factors $(\frac{1}{n}-1)(\frac{1}{n}-2)(\frac{1}{n}-3)\cdots$
has factorial growth which induces asymptotic convergence for this
representation. Because this is the well-known series of the confluent
hypergeometric function in (\ref{eq.hyp11}), we need not dwell into this further \cite[13.5.1]{AS}.

For small $x$, the Taylor series
\begin{equation}
\int \sin(x^n)dx
=\int \sum_{l=0}^\infty \frac{x^{n(2l+1)}}{(2l+1)!}(-)^ldx
=\sum_{l=0}^\infty \frac{(-)^l}{(2nl+n+1)} \frac{x^{2nl+n+1}}{(2l+1)!}
\label{eq.taylxn}
\end{equation}
and companion series for $\cos(x^n)$ converge quickly and are the standard target for Chebyshev approximations in that
area of numerical representation of special functions.

\begin{alg} \cite[p 63]{Erdelyi}
If the coefficients $b$ define a power series
\begin{equation}
f = \sum_{s\ge 0} b_s z^s,
\end{equation}
then the coefficients of the Neumann expansion
\begin{equation}
z^\nu f(z) =\sum_{s\ge 0}a_s J_{\nu+s}(z)
\end{equation}
are
\begin{equation}
a_s=(\nu+s)\sum_{m=0}^{\lfloor s/2\rfloor}
2^{\nu+s-2m}\frac{\Gamma(\nu+s-m)}{m!}b_{s-2m}
.
\end{equation}
\end{alg}
Application of this algorithm with $\nu=1+1/n$, $z=x^n$ and $n=2$ to Eq. (\ref{eq.taylxn})
yields the well established
\cite[8.515]{GR}\cite{BastardoAML18}
\begin{equation}
\int \sin(x^2)dx
=
\frac{\sqrt{2\pi}}{2}\sum_{s\ge 0}J_{2s+3/2}(x^2).
\end{equation}
For $n=3$ to $5$, the 
Neumann series
\begin{equation}
\int \sin(x^n)dx = d_n \sum_{s\ge 0} \xi_{n,s} J_{2s+\nu}(x^n)
\label{eq.Jn}
\end{equation}
are detailed in Tables \ref{tab.Jn3}--\ref{tab.Jn5} up to order $s=10$.
\begin{table}
\caption{Table of Neumann series coefficients (\ref{eq.Jn}) for $n=3$
with $d_3=2^{1/3}\pi \sqrt{3}/\Gamma(2/3) \approx 5.062876576879227$}
 \begin{tabular}{l|l}
$s$ & $\xi_{3,s}$\\
\hline
0 & 4/27 $\approx$ .14814814814814814814814815 \\ 
1 & 136/729 $\approx$ .18655692729766803840877915 \\ 
2 & 1120/6561 $\approx$ .17070568510897729004724889 \\ 
3 & 31912/177147 $\approx$ .18014417404754243650753329 \\ 
4 & 2491372/14348907 $\approx$ .17362799828586247022159946 \\ 
5 & 23052640/129140163 $\approx$ .17850867975131795365629204 \\ 
6 & 608995904/3486784401 $\approx$ .17465831951793224739736353 \\ 
7 & 16739336224/94143178827 $\approx$ .17780721272181224941292368 \\ 
8 & 148413541000/847288609443 $\approx$ .17516291302153317339294102 \\ 
9 & 12177019210000/68630377364883 $\approx$ .17742899977453386676382033 \\ 
10 & 325122763152640/1853020188851841 $\approx$ .17545559681899144476569162 \\ 
\end{tabular}
\label{tab.Jn3}
\end{table}
\begin{table}
\caption{Table of Neumann series coefficients (\ref{eq.Jn}) for $n=4$
with $d_4=2^{3/4}\pi/\Gamma(3/4) \approx 4.31160109908185588987751$}
 \begin{tabular}{l|l}
$s$ & $\xi_{4,s}$\\
\hline
0 & 1/8 $\approx$ .12500000000000000000000000 \\ 
1 & 11/64 $\approx$ .17187500000000000000000000 \\ 
2 & 159/1024 $\approx$ .15527343750000000000000000 \\ 
3 & 1347/8192 $\approx$ .16442871093750000000000000 \\ 
4 & 41531/262144 $\approx$ .15842819213867187500000000 \\ 
5 & 341309/2097152 $\approx$ .16274881362915039062500000 \\ 
6 & 5350155/33554432 $\approx$ .15944704413414001464843750 \\ 
7 & 43506995/268435456 $\approx$ .16207618638873100280761719 \\ 
8 & 2747379155/17179869184 $\approx$ .15991851425496861338615417 \\ 
9 & 22228087705/137438953472 $\approx$ .16173062398593174293637276 \\ 
10 & 352241109985/2199023255552 $\approx$ .16018071163898639497347176 \\ 
\end{tabular}
\label{tab.Jn4}
\end{table}
\begin{table}
\caption{Table of Neumann series coefficients (\ref{eq.Jn}) for $n=5$
with $d_5=2^{1/5}\pi/[\sin(\pi/5)\Gamma(4/5)] \approx 5.27349462002150700879064$}
 \begin{tabular}{l|l}
$s$ & $\xi_{5,s}$\\
\hline
0 & 2/25 $\approx$ .080000000000000000000000000 \\ 
1 & 72/625 $\approx$ .11520000000000000000000000 \\ 
2 & 8118/78125 $\approx$ .10391040000000000000000000 \\ 
3 & 214544/1953125 $\approx$ .10984652800000000000000000 \\ 
4 & 5179592/48828125 $\approx$ .10607804416000000000000000 \\ 
5 & 663616816/6103515625 $\approx$ .10872697913344000000000000 \\ 
6 & 16287458776/152587890625 $\approx$ .10674148983439360000000000 \\ 
7 & 2065603429328/19073486328125 $\approx$ .10829710907555184640000000 \\ 
8 & 51039744111764/476837158203125 $\approx$ .10703810144347409612800000 \\ 
9 & 1288445436120032/11920928955078125 $\approx$ .10808263692999989395456000 \\ 
10 & 159738815063405788/1490116119384765625 $\approx$ .10719890415611250403704832 \\ 
\end{tabular}
\label{tab.Jn5}
\end{table}
\begin{rem}
According to (\ref{eq.compl}),
$\lim_{s\to \infty} d_n\xi_{n,s} = \frac{2}{n}\Gamma(\frac{1}{n})\sin(\frac{\pi}{2n})
$.
The observation that $\xi_{n,s}$ stays practically constant as $s$ grows tells us
that this expansion converges poorly for medium and large upper limits of the integral.
\end{rem}

\section{Product with Powers, Incomplete Gamma Functions}\label{sec.powp}

Introducing a power of the free variable in the integrand leads from (\ref{eq.taylxn}) to
the more general
\begin{multline}
\int x^m \sin(x^n)dx
=\int \sum_{l=0}^\infty \frac{x^{m+2nl+n}}{(2l+1)!}(-)^ldx
=\sum_{l=0}^\infty \frac{(-)^l}{(m+2nl+n+1)} \frac{x^{m+2nl+n+1}}{(2l+1)!}
\\
= \frac{x^{m+n+1}}{m+n+1}
\,_1F_2\left(
\begin{array}{c}
\frac{1}{2}+\frac{m+1}{2n}
\\
\frac{3}{2}+\frac{m+1}{2n},\frac{3}{2}
\end{array}
\mid -\frac{x^{2n}}{4}
\right)
.
\end{multline}
This is the imaginary part of
the following confluent hypergeometric function \cite{Muller}:
\begin{multline}
\int x^m \exp(ix^n)dx
=\int \sum_{l=0}^\infty \frac{i^l x^{m+nl}}{l!} dx
=\sum_{l=0}^\infty \frac{i^l}{(m+nl+1)} \frac{x^{m+nl+1}}{l!}
\\
= \frac{x^{m+1}}{m+1}
\,_1F_1\left(
\begin{array}{c}
\frac{m+1}{n}
\\
1+\frac{m+1}{n}
\end{array}
\mid ix^n
\right)
.
\label{eq.hyp11}
\end{multline}
These underivatives are illustrated in Figure \ref{fig.Jn22to4S} and \ref{fig.Jn22to4C}.

\begin{rem}
By the rules $(a)_{2r}=(a)_r(\frac{1+a}{2})_r4^r$ and
$(a)_{2r+1}=a(\frac{1+a}{2})_r(1+\frac{a}{2})_r4^r$ for the Pochhammer Symbol \cite{SlaterHyp},
the even and odd parts of a generalized hypergeometric function $_rF_s(z)$ are \cite{ExtonJCAM79}
\begin{equation}
\sum_{n=0,2,4,6,\ldots} \frac{(a_1)_n (a_2)_n\cdots (a_r)_n}
{(b_1)_n (b_2)_n\cdots (b_s)_n}
\frac{z^n}{n!}
=\,_{2r}F_{2s+1}\left(\begin{array}{c}
\frac{a_1}{2},
\frac{1+a_1}{2},
\frac{a_2}{2},
\frac{1+a_2}{2},
\cdots
\frac{a_r}{2},
\frac{1+a_r}{2}\\
\frac{b_1}{2},
\frac{1+b_1}{2},
\frac{b_2}{2},
\frac{1+b_2}{2},
\cdots
\frac{b_r}{2},
\frac{1+b_r}{2},\frac{1}{2}\\
\end{array}\mid \frac{z^2}{4^{s-r+1}}\right)
\end{equation}
and
\begin{multline}
\sum_{n=1,3,5,7,\ldots} \frac{(a_1)_n (a_2)_n\cdots (a_r)_n}
{(b_1)_n (b_2)_n\cdots (b_s)_n}
\frac{z^n}{n!}
=
\\
\frac{\prod_{j=1}^r a_j}{\prod_{j=1}^s b_j}
z\,
_{2r}F_{2s+1}\left(\begin{array}{c}
\frac{1+a_1}{2},
1+\frac{a_1}{2},
\frac{1+a_2}{2},
1+\frac{a_2}{2},
\cdots
\frac{1+a_r}{2},
1+\frac{a_r}{2}\\
\frac{1+b_1}{2},
1+\frac{b_1}{2},
\frac{1+b_2}{2},
1+\frac{b_2}{2},
\cdots
\frac{1+b_r}{2},
1+\frac{b_r}{2},\frac{3}{2}\\
\end{array}\mid \frac{z^2}{4^{s-r+1}}\right),
\end{multline}
where parentheses with subindex denote Pochhammer's symbol,
\begin{equation}
(a)_n\equiv a(a+1)(a+2)\cdots (a+n-1)\equiv \Gamma(a+n)/\Gamma(a)
.
\end{equation}
Extraction of real and imaginary parts of hypergeometric functions with purely imaginary
argument $z$ is an application of these formulas.
\end{rem}
\begin{rem}\label{rem.reduO}
The increase in amplitude observed in these plots for large $x$
if $\phi$
has low degree $n$ does not hamper the numerical treatment.
The partial integration
\begin{equation}
in\int x^m e^{ix^n}dx = x^{m-n+1}e^{ix^n}-(m-n+1)\int x^{m-n}e^{ix^n}dx
\end{equation}
decrements the exponent $m$ in the integrand---closely related to a contiguous
relation of $_1F_1$ in (\ref{eq.hyp11}) \cite[13.4.7]{AS}.
Combined with (\ref{eq.triv}), the task that remains is to handle the cases $m<n$.
More generally speaking, if the degree of the polynomial $p(x)$ in the integrals of the
form $\int p(x)e^{i\phi(x)}dx$ is at least the degree of the polynomial $\phi(x)$ minus one,
a polynomial division of $p$ through the first derivative of $\phi$, $p(x)=\bar p(x)\phi'(x)+r(x)$,
followed by partial integration of the term $\int \bar p(x)\phi'(x)e^{i\phi(x)}dx$ reduces
the task to the format $\int \bar p'(x)e^{i\phi(x)}dx$ and $\int r(x)e^{i\phi(x)}dx$
where polynomial degrees are $\deg \bar p'<\deg \phi'$
and $\deg r < \deg \phi'$---this reduction executed
recursively if needed.
\end{rem}
The leading term in a Laurent expansion is \cite[13.5.1]{AS}
\begin{equation}
\,_1F_1\left(
\begin{array}{c}
\frac{m+1}{n}
\\
1+\frac{m+1}{n}
\end{array}
\mid ix^n
\right)
\to
\frac{m+1}{n}\Gamma(\frac{m+1}{n})
e^{i\pi (m+1)/(2n)}
x^{-(m+1)}
,
\end{equation}
which is inserted into (\ref{eq.hyp11}),
\begin{equation}
\int_0^\infty x^m\exp(ix^n)dx =
\frac{1}{n}\Gamma(\frac{m+1}{n})
e^{i\pi (m+1)/(2n)}
,
\label{eq.compl}
\end{equation}
and by the scaling substitution $xc^{1/n}\to y$,
\begin{equation}
\int_0^\infty x^m\exp(icx^n)dx
=
\frac{1}{nc^{(1+m)/n}}\Gamma(\frac{m+1}{n})
e^{i\pi (m+1)/(2n)}
 =
\frac{1}{n}(\frac{i}{c})^{(1+m)/n}\Gamma(\frac{m+1}{n})
.
\label{eq.complwc}
\end{equation}
\begin{rem}
Particular imaginary parts of these limits are 
\begin{equation}
\int_0^\infty \sin(x^n)dx = \frac{1}{n}\Gamma(1/n)\sin\frac{\pi}{2n}
\approx \begin{cases}
0.626657068657750125603941,&n=2; \\
0.446489755784624605609282,&n=3; \\
0.346865211023809496042035,&n=4. \\
\end{cases}
\end{equation}
\end{rem}
\begin{rem}
A simple pattern emerges for the $k$-fold repeated integrals \cite{DriverETNA25},
\begin{multline}
\int\cdots\int x^m \exp(ix^n)dx
=\sum_{l=0}^\infty \frac{i^l}{(m+nl+1)(m+nl+2)\cdots(m+nl+k)} \frac{x^{m+nl+k}}{l!}
\\
=\frac{x^{m+k}}{(m+1)(m+2)\cdots(m+k)}\,_kF_k\left(\begin{array}{cccc}
\frac{m+1}{n}, & \frac{m+2}{n}&\cdots & \frac{m+k}{n}
\\
1+\frac{m+1}{n}, & 1+\frac{m+2}{n} & \cdots & 1+\frac{m+k}{n}
\end{array}\mid ix^n\right)
.
\end{multline}
\end{rem}

Kummer's transformation of (\ref{eq.hyp11}) is \cite[13.1.27]{AS}
\begin{equation}
\int x^m \exp(ix^n)dx
\equiv V_{n,m}(x)e^{ix^n},
\label{eq.Vdef}
\end{equation}
with
\begin{equation}
V_{n,m}(x) = 
\frac{x^{m+1}}{m+1}
\,_1F_1\left(
\begin{array}{c}
1
\\
1+\frac{m+1}{n}
\end{array}
\mid -ix^n
\right)
.
\label{eq.Vnm}
\end{equation}
Real and imaginary part of (\ref{eq.Vdef}) are
\begin{equation}
\int x^m \cos(x^n)dx=
\Re V_{n,m}(x)\cos(x^n)
-\Im V_{n,m}(x)\sin(x^n);
\end{equation}
\begin{equation}
\int x^m \sin(x^n)dx=
\Re V_{n,m}(x)\sin(x^n)
+\Im V_{n,m}(x)\cos(x^n).
\end{equation}

\begin{rem}\label{rem.repPart}
The same result could have been derived by repeated
partial integration, the initial two steps being
\begin{multline}
\int x^m e^{ix^n}dx = \frac{1}{m+1}x^{m+1}e^{ix^n}-\frac{in}{m+1}\int x^{m+n}e^{ix^n}dx
\\
= \frac{1}{m+1}x^{m+1}e^{ix^n}
+\frac{-in}{(m+1)(m+n+1)}x^{m+n+1} e^{ix^n}
+\frac{(-in)^2}{(m+1)(m+n+1)}\int x^{m+2n}e^{ix^n}dx
.
\end{multline}
\end{rem}

Gathering each second term of the hypergeometric series of (\ref{eq.Vnm}),
\begin{multline}
\Im V_{n,m}(x)
=
-\frac{x^{m+1}}{m+1}\sum_{k=0}^{\infty} \frac{(-)^k(x^n)^{2k+1}}{(1+\frac{m+1}{n})_{2k+1} (2k+1)!}
\\
=
-\frac{x^{m+n+1}}{(m+1)(1+\frac{m+1}{n})}
\,_1F_2\left(
\begin{array}{c}
1 \\
1+\frac{m+1}{2n}, \frac{3}{2}+\frac{m+1}{2n}
\end{array}
\mid -\frac{x^{2n}}{4}
\right)
.
\end{multline}

A corresponding calculation of the other half of the terms yields
\begin{equation}
\Re V_{n,m}(x)
=
\frac{1}{m+1}x^{m+1}
\,_1F_2\left(\begin{array}{c}
1\\
\frac{1}{2}+\frac{m+1}{2n}, 1+\frac{m+1}{2n}
\end{array}
\mid -\frac{x^{2n}}{4}\right)
\label{eq.VWfac}
.
\end{equation}
\begin{rem}\label{rem.anti}
One might hope that splitting off the factor $\exp(ix^n)$ in
(\ref{eq.Vdef}) removes the major oscillations and $V_{n,m}$ contains
less wiggles. Figure \ref{fig.VWfac} illustrates
that this
target is missed, and Section \ref{sec.ancy} puts this into a wider perspective.
\begin{figure}
\includegraphics[width=0.8\textwidth]{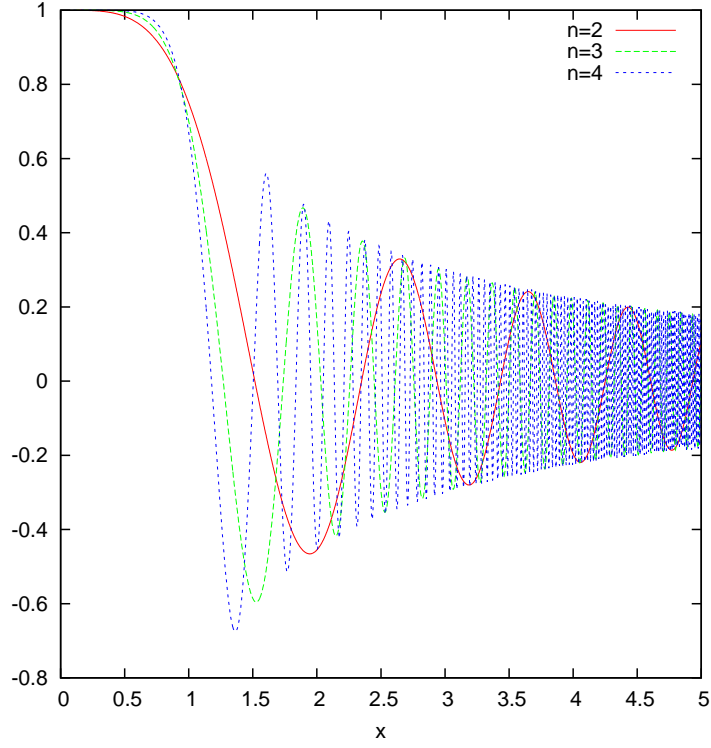}
\caption{The value of
$_1F_2(1;\frac{1}{2}+\frac{1}{2n},1+\frac{1}{2n};-\frac{x^{2n}}{4})$,
the main constituent of $\Re V_{n,0}(x)$ in (\ref{eq.VWfac}),
for
three different degrees $n$ of the phase polynomial.}
\label{fig.VWfac}
\end{figure}
\end{rem}

\begin{rem}
The continued fraction of (\ref{eq.Vnm}) is \cite{Cuyt,Muller,Jones}
\begin{multline}
_1F_1\left(\begin{array}{c}1 \\ 1+a\end{array}\mid x\right)
=
1+\frac{x}{\frac{1+a}{0!}+} \,
\frac{x}{ -\frac{0!(2+a)}{(1+a)}+}
\frac{x}{ -\frac{(1+a)(3+a)}{1!}+}
\frac{ x}{ \frac{1!(4+a)}{(1+a)(2+a)} + }
\\
\frac{x}{ \frac{(1+a)(2+a)(5+a)}{2!}+}
\frac{x}{ -\frac{2!(6+a)}{(1+a)(2+a)(3+a)} +}
\frac{x}{ -\frac{(1+a)(2+a)(3+a)(7+a)}{3!}+}
\frac{x}{ \frac{3!(8+a)}{(1+a)(2+a)(3+a)(4+a)} +}
\\
\frac{x }{ \frac{(1+a)(2+a)(3+a)(4+a)(9+a)}{4!}+}
\frac{ x}{-\cdots}
\end{multline}

Fields and Wimp and Luke provide a series expansion in terms of Bessel functions
\cite[(2.14)]{Fields}\cite{Luke1959}.
\end{rem}

Eq. (\ref{eq.hyp11}) transforms into an Incomplete Gamma Function
\cite[13.6.10]{AS}\cite[2.632]{GR}\cite{AmoreEPLett71,GautschiMCom31,GautschiCiteseer,Chaudhry2002,Allasia,Barakat}:
\begin{equation}
\int x^m \exp(ix^n)dx
= 
\frac{1}{n}
i^{(m+1)/n}
\gamma(\frac{m+1}{n},-ix^n)
,
\end{equation}
which is sometimes phrased as
\cite[6.5.3]{AS}
\begin{equation}
\gamma(a,z)=\Gamma(a)-\Gamma(a,z)
.
\label{eq.ggg}
\end{equation}

This association opens a wealth of literature to numerical evaluation.
A particularly nicely converging continued fraction representation is \cite[12.6.34]{Cuyt}\cite{TemmeISNM119}
\begin{equation}
\Gamma(s,z) = \frac{z^se^{-z}}{1+z-s+}\,
\frac{s-1}{3+z-s+} \,
\frac{2(s-2)}{5+z-s+} \,
\frac{3(s-3)}{7+z-s+} \,
\frac{4(s-4)}{9+z-s+} \,\cdots
\end{equation}
\begin{rem}
This can apparently be re-binned in terms of 
associated Laguerre polynomials as \cite{Crandall}\cite[(41)]{TricomiMA31}
\begin{equation}
\Gamma(s,x)=x^s e^{-x}\sum_{n\ge 0}\frac{(1-s)_n}{(n+1)!}\,\frac{1}{L_n^{-s}(-x) L_{n+1}^{-s}(-x)}.
\end{equation}
\end{rem}

\section{Series Reversion}\label{sec.rev}
The final chapters of the manuscript are concerned with
phases $\phi(x)$ which are not integer powers but full-fledged polynomials of $x$.
\subsection{Expansion Around the Origin}\label{sec.serO}
If the leading (lower) degree of the phase polynomial is 1 (so $\alpha_0=0$
after use of (\ref{eq.sina0}) and (\ref{eq.cosa0}) and $\alpha_1\neq 0$),
a substitution
\begin{equation}
y\equiv \sum_{j=1}^l\alpha_j x^j
\label{eq.yofr}
\end{equation}
in the integrand is available, associated with a series
reversion of the polynomial,
\begin{eqnarray}
x &\equiv& \sum_{j=1}^\infty \beta_j y^j, \label{eq.xofy}\\
dx/dy &=& \sum_{j=1}^\infty j\beta_j y^{j-1}.
\end{eqnarray}
The idea is that preserving the sine/cosine/exponential in the integrand
allows to deal efficiently with the oscillations in the spirit of Filon quadratures.

\begin{alg}
The intermediate task is to determine the (infinite) set of the 
$\beta_j$ from the (finite set of) $\alpha_j$.
The $\beta_j$ are tabulated as a function of the $\alpha_j$
by Abramowitz and Stegun \cite[3.6.25]{AS} for $j\le 7$,
and by Orstrand \cite{OrstrandPM37} for $j\le 13$. In the general case, equalize
the coefficients of equal powers of $y$ in (\ref{eq.yofr}), which is
\begin{multline}
y=\sum_{j=1}^l\alpha_j\left(\sum_{k=1}^\infty \beta_ky^k\right)^j
\\
=
\sum_{j=1}^l\alpha_j \sum_{t=j}^\infty y^t
\sum_{\begin{array}{c} t_1+2t_2+\cdots +mt_m=t\\ t_1+t_2+\cdots +t_m=j\\ t_1, t_2,\ldots ,t_m\ge 0\end{array}}
\left(\begin{array}{c}j\\ t_1\, t_2\, \ldots t_m\end{array}\right)
\beta_1^{t_1}
\beta_2^{t_2}
\cdots
\beta_m^{t_m}
\end{multline}
involving the multinomial coefficients in parentheses \cite{RaneyTAMS94,OrstrandPM37,ChangAMC23,FielderMCom14,WardRCM54,Dominiciarxiv05}.
The maximum index $m$ that contributes with $t_m>0$ to the
sums $\sum_1^m it_i=t$ and $\sum_1^m t_i=j$ is
$m=t-j+1$.
The number of multinomial terms can be read off entry A008284 in the
Online Encyclopedia of Integer Sequences (OEIS) \cite{EIS}.
Their evaluation may be avoided by employing a convolutional recurrence \cite[0.314]{GR}.

In the order $O(y^1)$,
\begin{equation}
1=\alpha_1\beta_1.
\end{equation}
In the orders $O(y^t)$, $t>1$, $\beta_t$
is derived by isolating $\beta_t$ in the term with $j=1$ in
\begin{multline}
0
=
\sum_{j=1}^{\min(t,l)}\alpha_j
\sum_{\begin{array}{c} t_1+2t_2+\cdots +mt_m=t\\ t_1+t_2+\cdots +t_m=j\\ t_1, t_2,\ldots, t_m\ge 0\end{array}}
\left(\begin{array}{c}j\\ t_1\, t_2\, \ldots t_m\end{array}\right)
\beta_1^{t_1}
\beta_2^{t_2}
\cdots
\beta_m^{t_m}
\\
=
\alpha_1
\sum_{\begin{array}{c} t_1+2t_2+\cdots +mt_m=t\\ t_1+t_2+\cdots +t_m=1\\ t_1, t_2,\ldots, t_m\ge 0\end{array}}
\left(\begin{array}{c}1\\ 0\, 0\, \ldots 0\, 1\end{array}\right)
\beta_1^0
\beta_2^0
\cdots
\beta_t^1
\\
+
\sum_{j=2}^{\min(t,l)}\alpha_j
\sum_{\begin{array}{c} t_1+2t_2+\cdots +mt_m=t\\ t_1+t_2+\cdots +t_m=j\\ t_1, t_2,\ldots, t_m\ge 0\end{array}}
\left(\begin{array}{c}j\\ t_1\, t_2\, \ldots t_m\end{array}\right)
\beta_1^{t_1}
\beta_2^{t_2}
\cdots
\beta_m^{t_m}
, t>1.
\end{multline}
which yields the recursive algorithm
\begin{equation}
\alpha_1 \beta_t
=
-\sum_{j=2}^{\min(t,l)}\alpha_j
\sum_{\begin{array}{c} t_1+2t_2+\cdots +mt_m=t\\ t_1+t_2+\cdots +t_m=j\\ t_1, t_2,\ldots, t_m\ge 0\end{array}}
\left(\begin{array}{c}j\\ t_1\, t_2\, \ldots t_m\end{array}\right)
\beta_1^{t_1}
\beta_2^{t_2}
\cdots
\beta_m^{t_m}
, t>1.
\end{equation}
In practise, the problem is solved by a perturbative expansion,
starting from the estimate $x^{(1)}=y/\alpha_1$
and inserting iteratively previous estimates into the constituent equation:
$x^{(s)}=y-\sum_{j=2}^l\alpha_j (x^{(s-1)})^j$ \cite{BrentJACM25}.
\end{alg}

For a set of simple polynomials, the $\beta$ coefficients are gathered in Appendix \ref{sec.revSimpl}.

The application of the series reversion linearizes the phase argument in the
trigonometric function.
\begin{exa}
\begin{multline}
\int x\exp(i\sum_{j=1}^l \alpha_j x^j) dx
=
\int (\sum_{k=1}^\infty k\beta_k y^{k-1})(\sum_{j=1}^\infty \beta_j y^j)\exp(iy) dy,
\\
=
\sum_{t=1}^\infty \gamma_t \int y^t\exp(iy) dy,
\label{eq.revsd}
\end{multline}
where the series coefficients $\gamma$ are derived from convolution
type products
\begin{equation}
\gamma_t\equiv \sum_{k=1}^t k\beta_k \beta_{t-k+1}
.
\end{equation}
\end{exa}
Approximation of the representation of the series 
by truncation at some maximum
index of the $\beta$ reduces the problem
for each term in (\ref{eq.revsd}) to the closed forms (\ref{eq.xnsinax}).

\begin{rem} \label{rem.sqr}
If the polynomial in the phase argument has degree $l=2$, the linear term
in the polynomial $\alpha_0+\alpha_1x+\alpha_2x^2$ may be eliminated by the substitution
$z=x+\alpha_{l-1}/(l\alpha_l)$. The remaining absolute term can also be
removed as mentioned in Remark \ref{rem.o0}, so
$\int \cos(\alpha_0+\alpha_1x+\alpha_2x^2)dx$
and
$\int \sin(\alpha_0+\alpha_1x+\alpha_2x^2)dx$
have closed form expressions in terms of Fresnel Integrals
\cite[2.549.3, 2.549.4]{GR}.
\end{rem}

\subsection{Local Series}
The main reason \emph{not} to use the approach of Section \ref{sec.serO}
in numerical application
is that the series reversion is only unique/stable up to the argument
of the first maximum of the phase polynomial $y(x)$. 
The administrative problem is to monitor
branch cuts of the inverse $x(y)$ and to monitor the stability of the 
results as a function of the truncation order of the $\beta$.
The quick growth of the sequences in Appendix \ref{sec.revSimpl} demonstrates the problems
of convergence, and for small upper limits within the radius of convergence,
the competing Taylor series are simpler to handle.
One improvement is to split the interval of integration
and to revert Taylor series of the phase polynomial anchored
in the middle of the sub-intervals.

The power series expansion of the integrand near some $x$
is found with the binomial expansion
\begin{multline}
\exp(i(x+\epsilon)^n) =
\exp\left(ix^n+i\sum_{l=1}^n\binom{n}{l}\epsilon^lx^{n-l}\right)
\\
=
\exp(ix^n)\exp\left(i\sum_{l=1}^n\binom{n}{l}\epsilon^lx^{n-l}\right)
.
\label{eq.bino}
\end{multline}
Assuming $\epsilon$ is small, reversion of the
polynomial in
the
argument on the right hand side
generates the power series. We illustrate this procedure in the
next three subsections.

\subsubsection{Squares}
For quadratic chirp polynomials, $n=2$, 
\begin{equation}
\exp\left(i\sum_{l=1}^n\binom{n}{l}\epsilon^lx^{n-l}\right)
=\exp(2ix\epsilon+i\epsilon^2) = \exp(2ixy),
\end{equation}
where 
\begin{equation}
y\equiv \epsilon+\frac{1}{2x}\epsilon^2
\end{equation}
is substituted in the integral.
Series reversion
\begin{equation}
\epsilon = y
-\frac{1}{2x}y^2
+\frac{2}{(2x)^2}y^3
-\frac{5}{(2x)^3}y^4
+\frac{14}{(2x)^4}y^5-\ldots
\end{equation}
with coefficients from the first row of Table \ref{tab.revs}
leads to the derivative
\begin{equation}
d\epsilon/dy = 
1
-\frac{2}{2x}y
+\frac{2\cdot 3}{(2x)^2}y^2
-\frac{4\cdot 5}{(2x)^3}y^3
+\frac{5\cdot 14}{(2x)^4}y^4-\ldots
=\sum_{l\ge 0} \frac{\binom{2l}{l}}{(2x)^l} (-y)^l
.
\end{equation}
Integration over small $\epsilon$-environments near some $x$ are
then evaluated with
\begin{equation}
\int \exp(2ix\epsilon+i\epsilon^2)d\epsilon
=\sum_{l\ge 0} \frac{\binom{2l}{l}}{(2x)^l}
\int
(-y)^l\exp(2ixy) dy
\label{eq.coseps2}
\end{equation}
--- eventually with (\ref{eq.xnsinax})---.
In practise, these quadratic polynomials would rather be evaluated with 
the closed form procedure of Remark \ref{rem.sqr}.

\subsubsection{Cubes}
For cubic chirp polynomials, application of (\ref{eq.bino}) reads
\begin{equation}
\exp(i(x+\epsilon)^3)=\exp(ix^3)\exp(3i\epsilon x^2+3i\epsilon^2x+i\epsilon^3)
\label{eq.cubs}
\end{equation}
with one factor
\begin{equation}
\exp(3i\epsilon x^2+3i\epsilon^2x+i\epsilon^3)
=
\exp(3ix^2(\epsilon +\frac{1}{x}\epsilon^2+\frac{1}{3x^2}\epsilon^3))
=
\exp(3ix^2y)
\end{equation}
after substituting
\begin{equation}
y=\epsilon+\frac{\surd 3}{\surd 3 x}\epsilon^2+\frac{1}{3x^2}\epsilon^3.
\label{eq.cubsy}
\end{equation}
This reverts to \cite[A218540]{EIS}
\begin{multline}
\epsilon = y
-\frac{\surd 3}{\surd 3 x}y^2
+\frac{5}{(\surd 3 x)^2}y^3
-\frac{10 \surd 3}{(\surd 3 x)^3}y^4
+\frac{66}{(\surd 3 x)^4}y^5
-\frac{154 \surd 3}{(\surd 3 x)^5}y^6
+\frac{1122}{(\surd 3 x)^6}y^7
-\frac{2805 \surd 3}{(\surd 3 x)^7}y^8
+\ldots
\\
= \sum_{j\ge 0} \frac{\kappa_j}{(\surd 3 x)^{2j}}y^{2j+1}
-
\sum_{j\ge 0} \frac{\lambda_j \surd 3}{(\surd 3 x)^{2j+1}}y^{2j+2}
.
\label{eq.cubsser}
\end{multline}
\begin{rem}
Alternatively, by canceling some square roots of 3,
\begin{equation}
\epsilon = y
-\frac{1}{x}y^2
+\frac{15}{(3 x)^2}y^3
-\frac{90}{(3x)^3}y^4
+\frac{594}{(3x)^4}y^5
-\frac{4158}{(3 x)^5}y^6
+\ldots
\end{equation}
with numerator coefficients tabulated in the
OEIS sequence A025748 \cite{EIS}.
\end{rem}
On the r.h.s., the coefficients $\kappa_0=1$, $\kappa_1=5$, $\kappa_2=66$, $\kappa_3=1122$ etc\@. 
established in the
numerators of odd powers of $y$ obey $n(2n+1)\kappa_n = 3(6n-1)(3n-2)\kappa_{n-1}$.
The recurrence is solved by
\begin{equation}
\kappa_n = 27^n (5/6)_n\, (1/3)_n/[n! (3/2)_n]
\end{equation}
with generating function
\begin{equation}
\sum_{n\ge 0}\kappa_n z^n=\,_2F_1\left(\begin{array}{c}5/6, 1/3\\ 3/2\end{array}\mid 27z\right)
.
\end{equation}
The coefficients $\lambda_0=1$, $\lambda_1=10$, $\lambda_2=154$, $\lambda_3=2805,\ldots$ in the
numerators of even powers of $y$ obey $(n+1)(2n+1)\lambda_n = 3(6n-1)(3n+1)\lambda_{n-1}$,
such that
\begin{equation}
\lambda_n = \frac{27^n (5/6)_n\, (4/3)_n}{(n+1)! (3/2)_n}
=\frac{1+3n}{1+n}\kappa_n
\end{equation}
with generating function
\begin{equation}
\sum_{n\ge 0}\lambda_n z^n = 
\,_3F_2\left(\begin{array}{c}5/6, 4/3, 1\\ 3/2, 2\end{array}\mid 27z\right)
.
\end{equation}
The substitution actually calls for the derivatives
\begin{equation}
d\epsilon/dy 
= \sum_{j\ge 0} \frac{(2j+1)\kappa_j}{(\surd 3 x)^{2j}}y^{2j}
-
\sum_{j\ge 0} \frac{(2j+2)\lambda_j \surd 3}{(\surd 3 x)^{2j+1}}y^{2j+1}
.
\label{eq.cubsdiff}
\end{equation}
So the cubic analog of (\ref{eq.coseps2}) is
\begin{equation}
\int
\exp(3i\epsilon x^2+3i\epsilon^2x+i\epsilon^3) d\epsilon
=
\int dy
\exp(3ix^2y) 
[\sum_{j\ge 0} \frac{(2j+1)\kappa_j}{(\surd 3 x)^{2j}}y^{2j}
-
\sum_{j\ge 0} \frac{(2j+2)\lambda_j \surd 3}{(\surd 3 x)^{2j+1}}y^{2j+1}]
\label{eq.cubsrev}
\end{equation}
for use in the right hand side of (\ref{eq.bino}).

\subsubsection{Quartics}
The quartic chirp with (\ref{eq.bino}) reads
\begin{equation}
\exp(i(x+\epsilon)^4)=\exp(ix^4)\exp(4i\epsilon x^3+6i\epsilon^2x^2+4i\epsilon^3x + i\epsilon^4)
\end{equation}
with small $y$ in
\begin{equation}
\exp(4i\epsilon x^3+6i\epsilon^2x^2+4i\epsilon^3x + i\epsilon^4)
=
\exp\left[4ix^3(\epsilon +\frac{3}{2x}\epsilon^2+\frac{1}{x^2}\epsilon^3+\frac{1}{4x^3}\epsilon^4)\right]
=
\exp(4ix^3y)
\end{equation}
substituting
\begin{equation}
y=\epsilon+\frac{3}{2x}\epsilon^2+\frac{1}{x^2}\epsilon^3+\frac{1}{4x^3}\epsilon^4.
\end{equation}
This reverts to
\begin{multline}
\epsilon = y
-\frac{3}{2^1x}y^2
+\frac{14}{2^2x^2}y^3
-\frac{77}{2^3x^3}y^4
+\frac{462}{2^4x^4}y^5
-\frac{2926}{2^5x^5}y^6
+\frac{19228}{2^6x^6}y^7
-\frac{129789}{2^7x^7}y^8
+\ldots
\\
= \sum_{j\ge 1} (-)^{j+1} \frac{\eta_j}{(2x)^{j-1}}y^j,
\end{multline}
where the integer sequence $\eta$ is built up by  \cite[A048779]{EIS}
\begin{equation}
j\eta_j+2(5-4j)\eta_{j-1}=0;\quad \sum_{j\ge 1}\eta_j z^j=\frac{1-\sqrt[4]{1-8z}}{2}.
\end{equation}
With derivatives
\begin{equation}
d\epsilon/dy 
= \sum_{j\ge 0} (-)^j\frac{(j+1)\eta_{j+1}}{(2x)^j}y^j
,
\end{equation}
the quartic variant of (\ref{eq.coseps2}) is
\begin{equation}
\int
\exp(4i\epsilon x^3+6i\epsilon^2x^2+4i\epsilon^3x + i\epsilon^4)
d\epsilon
=
\int dy
\exp(4ix^3y) \sum_{j\ge 0} \frac{(j+1)\eta_{j+1}}{(2x)^j}(-y)^j
.
\end{equation}

\section{Associated Linear Differential Equation}\label{sec.pract}
\subsection{Perturbative}

The iterated partial integration of Remark \ref{rem.repPart}
defines a generalization of (\ref{eq.Vdef}) to polynomials $p(x)$ and $\phi(x)$,
\begin{equation}
I_{p,\phi}(x) = \int p(x)e^{i\phi(x)}dx
=q_1(x)e^{i\phi(x)}-\int iq_1(x)\phi' e^{i\phi}dx
= e^{i\phi(x)}\sum_{k=0}^\infty q_{2k+1}(x),
\end{equation}
where
\begin{equation}
q_0(x)\equiv p(x), \quad q_1(x)\equiv \int q_0(x)dx,\quad
q_{2k+1}(x)\equiv -i\int q_{2k-1}(x)\phi'(x)dx,\quad k>0
\label{eq.qiter}
\end{equation}
construct a sequence of polynomials of increasing degree.
The sum $\sum q_{2k+1}$ equals $V_{n,m}$ if $p$ and $\phi$ are integer powers of $x$;
in that sense $q(x)$ constructs the Kummer transformation of $I_{p,\phi}(x)$.
The sum over all
$q_{2k+1}$ is
a perturbative construction of a solution $q(x)$ to
the first-order linear inhomogeneous differential equation
\begin{equation}
p(x) = iq(x)\phi'(x)+q'(x),
\label{eq.Levin}
\end{equation}
respectively integral equation
\begin{equation}
q(x) = \int_0 p(x)dx - \int_0 i q(x)\phi'(x)dx.
\end{equation}
It starts with the estimate that $iq\phi'$ on the right hand side vanishes, which
leads to the approximation that $q$ is the integral over $p$. This is plugged into
the first term of the right hand side, and solved by integration over $p-iq\phi'$
and so on. We are mainly interested in singular solutions $I_{p,\phi}$ with lower limit zero, so 
we insert lower limits of zero also in the definitions of the integrals in (\ref{eq.qiter}).
\begin{rem}
The solution to the linear differential equation (\ref{eq.Levin}) is the sum over
the set of solutions with inhomogeneous terms $p(x)$ which are powers of $x$. These are Bernoulli
differential equations, so the standard substitution
$f=q^{1-\deg p}$ reduces the problem to the differential equations
\begin{equation}
\frac{df}{dx} + i(1-\deg p)\phi'(x)f(x) = 1-\deg p.
\end{equation}
In that sense, the only form to be studied is the one with constant left hand side in (\ref{eq.Levin}).
This route has not been followed here, because the connection to the initial value problem
$\lim_{x\to 0}q\to 0$ is lost by changing the variable to $f$ which is an inverse power of $q$.
\end{rem}

\begin{figure}
\includegraphics[width=0.8\textwidth]{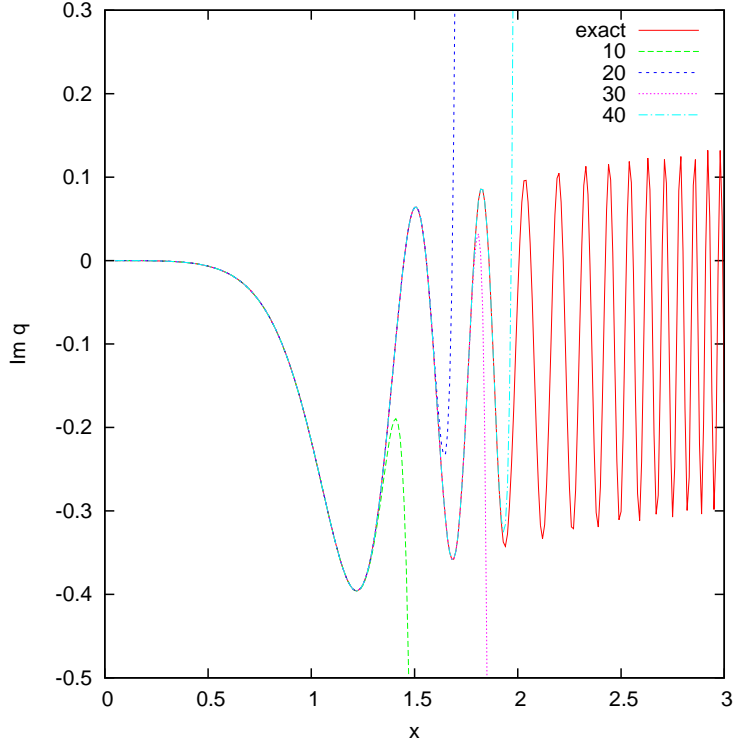}
\caption{The exact value of $\Im q_{p,\phi}$ as a function of upper limit
$x$
for the example $p=x^2$, $\phi=x+x^4$, and four polynomial approximations
that gather 10, 20, 30 or 40 orders of the expansion (\ref{eq.qiter}).
The approximations stay on the solution for approximately 1 to 3 oscillations
and exit abruptly.
}
\label{fig.ada}
\end{figure}

\begin{rem}
This is an analytic variant of Levin's collocation method of
casting
an approximate solution \cite{LevinMC38,ChungANM34}. It is geared to
a power-series solution (the degree of the polynomials $q$ increases with their index)
and in that sense
another
method to generate the Taylor series at $x=0$.
\end{rem}

\subsection{Asymptotics}
\begin{figure}
\includegraphics[width=0.8\textwidth]{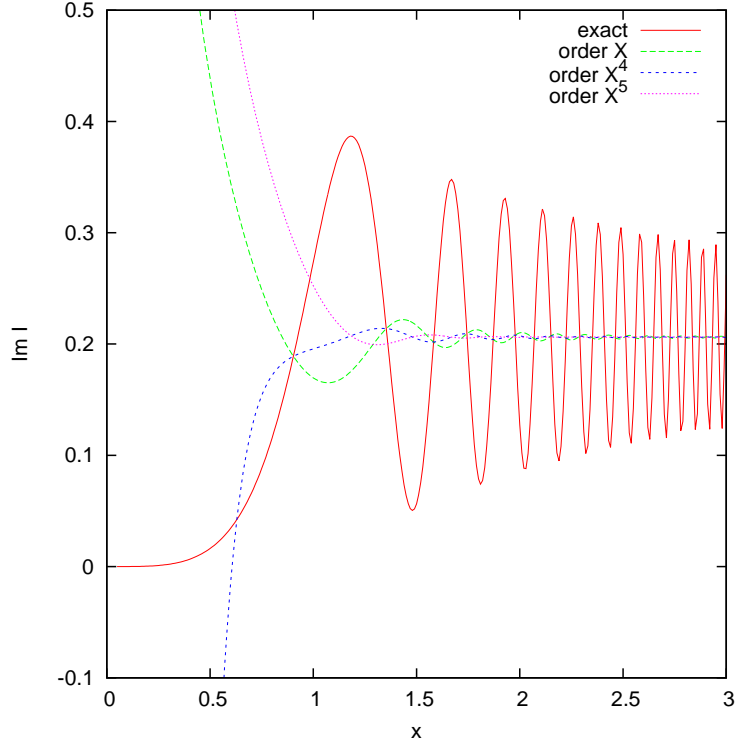}
\caption{The exact value of $\Im I_{p,\phi}$ as a function of upper limit
$x=1/X$
for the polynomials $p$ and $\phi$ of Example \ref{ex.asy}, plus three residual 
asymptotically horizontal curves
$\Im [I-e^{i\phi(x)}(-iX/4)]$,
$\Im[I_{p,\phi}-e^{i\phi(x)}(-iX/4+iX^4/16)]$
and
$\Im[I_{p,\phi}-e^{i\phi(x)}(-iX/4+iX^4/16-X^5/16)]$
after up to three leading components of the asymptotic
have been removed.
}
\label{fig.asy}
\end{figure}

The asymptotic behavior
\begin{equation}
q\sim \sum_{j=1}^\infty \frac{h_j}{x^j}=\sum_{j\ge 1}h_j X^j
\label{eq.asy}
\end{equation}
is found by substituting $X\equiv x^{-1}$ in the polynomials $p$ and $\phi'$,
replacing $dq/dx= -X^2 dq/dX$ in (\ref{eq.Levin}),
multiplication of the equation by $X^{\deg \phi'}$ and solving the differential
equation for $q(X)$ with a power series ansatz \cite{TurrittinFE12}, see Appendix \ref{app.asy}.
The leading term
of the series has order $X^{\deg \phi'-\deg p}$, where we
may assume that this series exists by enforcing $\deg\phi' > \deg p$
as discussed in Remark \ref{rem.reduO}.
\begin{exa}\label{ex.asy}
If $p=x^2$ and $\phi=x+x^4$, $\phi'=1+4x^3$, the differential equation
is $X^{-2}=iq(1+4X^{-3})-X^2q'(X)$, or $X=iq(X^3+4)-X^5q'(X)$,
which is solved by
\begin{equation}
q \sim -\frac{i}{4x}+\frac{i}{16x^4}-\frac{1}{16x^5}-\frac{i}{64x^7}+\cdots.
\end{equation}
The efficiency of this series to estimate the amplitude
of the oscillations of $I_{p,\phi}$ for large $x$
is illustrated in Figure \ref{fig.asy}.
\end{exa}

\subsection{Anticyclic}\label{sec.ancy}

So far, the expansion (\ref{eq.qiter})
has anchored the values of $I_{p,\phi}(x)$ at $x=0$,
and the asymptotic expansion (\ref{eq.asy})
has anchored them at $x=\infty$. The two solutions $q$
of different initial conditions miss each other by a gap $I_{p,\phi}(\infty)$,
\begin{equation}
e^{i\phi(x)}\sum_{k\ge 0}q_{2k+1} = I_{p,\phi}(\infty)+e^{i\phi(x)}\sum_{j\ge 1} h_jX^j.
\label{eq.Iinf}
\end{equation}
The solutions differ by a multiple of 
$I_{p,\phi}(\infty) e^{-i\phi}$, the solution to the homogeneous differential equation,
following 
the general theory of linear differential equations.
\begin{rem}
$I_{p,\phi}(\infty)$ might be called the ``complete'' generalized Fresnel Integral.

If we would start the perturbative accumulation of $q$ not with the estimate $p$ as
in (\ref{eq.qiter}) but with some constant $q_0=\epsilon$, the repeated integration and
multiplication with $\phi'$ in (\ref{eq.qiter}) would generate the term $-i\int \epsilon\phi'dx = -i\epsilon\phi$,
then
the term $-i\int (-i\epsilon\phi)\phi'dx = \epsilon(-i\phi)^2/2$, and eventually the
entire general solution $\epsilon e^{-i\phi}$ of the homogeneous differential equation.
The equivalent procedure applied to (\ref{eq.qiter}), splitting $q_{2k-1}$ at each iteration
into a multiple of $\phi$ and remainder, yields a representation $q=\sum_{s\ge 0}\bar q_s(x)(i\phi)^s/s!$,
where $\bar q_s$ are polynomials of degree less than $\deg \phi$.
\end{rem}

The effect on $q$ is that $\sum q_{2k+1}$ picks up an anti-resonant,
counter-cycling
phase $I_{p,\phi}(\infty) e^{-i\phi(x)}$ as $x$ increases, mentioned
in Remark \ref{rem.anti} and observed again in Figure \ref{fig.ada}.

The calculation of $I(\infty)$ is
accomplished
by Taylor expansion of all contributions
of $\phi$, setting apart the term of highest polynomial degree, using (\ref{eq.complwc}).
\begin{alg}
Let $\phi=\sum_{j=1}^l \alpha_j x^j$, $l=\deg \phi$. Obtain Taylor series coefficients
$t$ by letting
$p(x)e^{i\sum_{j=1}^{l-1} \alpha_j x^j}=\sum_{j\ge 0}t_jx^j$
---optionally with an intermediate expansion of the exponential factor
in complete exponential Bell polynomials \cite{KolbigJCAM69}.
Then
\begin{equation}
lI_{p,\phi}(\infty)=\sum_{j\ge 0} t_j
\Gamma(\frac{1+j}{l})
\left(\frac{i}{\alpha_l}\right)^{(1+j)/l}
\label{eq.Iinfgen}
.
\end{equation}
\end{alg}
Results of such calculation are put into Table \ref{tab.Iinf}. Cases
that can be reduced with Remark \ref{rem.reduO} or Remark \ref{rem.sqr}
or Eq.\ (\ref{eq.eodd})
have been excised.
\begin{table}
\caption{Examples of real and imaginary parts of
integrals $I_{p,\phi}$ at infinite upper limit.
}
 \begin{tabular}{l|l|r|r}
$p$ & $\phi$ & $\Re I_{p,\phi}(\infty)$ & $\Im I_{p,\phi}(\infty)$ \\
\hline
$1$ & $x+x^3$ & 0.41494101283606350 & 0.53411593027204143\\
$x$ & $x+x^3$ & -0.02016219157077424 & 0.30316529215034456\\
$1$ & $x-x^3$ & 1.11098231192364267 & -0.03858690498389667\\
$x$ & $x-x^3$ & 0.59330541726382226 & -0.22202080248217837\\
$1$ & $x^2+x^3$ & 0.54028350983057729 & 0.40844024533897794\\
$x$ & $x^2+x^3$ & 0.07962821700232091 & 0.26244025603540361\\
$1$ & $x^2-x^3$ & 1.25860675774543236 & -0.27494486339726367\\
$x$ & $x^2-x^3$ & 0.73476055145142588 & -0.44097066691842955\\
$1$ & $-x^2+x^3$ & 1.25860675774543236 & 0.27494486339726367\\
$x$ & $-x^2+x^3$ & 0.73476055145142588 & 0.44097066691842955\\
$1$ & $-x^2-x^3$ & 0.54028350983057729 & -0.40844024533897794\\
$x$ & $-x^2-x^3$ & 0.07962821700232091 & -0.26244025603540361\\
$1$ & $x+x^4$ & 0.50249679573307246 & 0.52601826686302930\\
$x$ & $x+x^4$ & 0.05824661765885730 & 0.32482243282841976\\
$x^2$ & $x+x^4$ & -0.07182498024997954 & 0.20623764810943091\\
$1$ & $x-x^4$ & 1.04843033186615915 & 0.09812564203676180\\
$x$ & $x-x^4$ & 0.55603788569756914 & -0.06487860797927502\\
$x^2$ & $x-x^4$ & 0.37604412360259629 & -0.15808966880512279\\
$1$ & $x^2+x^4$ & 0.60419228699943584 & 0.38964387635793124\\
$x^2$ & $x^2+x^4$ & -0.00678784821525451 & 0.18281575355157232\\
$1$ & $x+x^5$ & 0.56055180704649184 & 0.51785573267722635\\
$x$ & $x+x^5$ & 0.11195277049251853 & 0.33139696436598476\\
$x^2$ & $x+x^5$ & -0.02225695591532138 & 0.22710191687546598\\
$x^3$ & $x+x^5$ & -0.08199028168861113 & 0.15444851391776769\\
$1$ & $2x+x^3$ & 0.18186049037842599 & 0.46158113697480635\\
$x$ & $2x+x^3$ & -0.10324653677125750 & 0.16619667306199227\\
$1$ & $x^3+x^4$ & 0.66346454706590291 & 0.32615196348829038\\
$x$ & $x^3+x^4$ & 0.18566011304028281 & 0.23801625874900127\\
$x^2$ & $x^3+x^4$ & 0.03767486947451176 & 0.17284523624119006\\
$1$ & $x^2+x^5$ & 0.64955004591503774 & 0.37714240719158655\\
$x$ & $x^2+x^5$ & 0.17628609937611991 & 0.26821665359731032\\
$x^2$ & $x^2+x^5$ & 0.02993102573558424 & 0.19494632498226022\\
$x^3$ & $x^2+x^5$ & -0.03639133933981598 & 0.13909435860476519\\
$1$ & $x+x^6$ & 0.60171369622318495 & 0.51091411953631571\\
$x$ & $x+x^6$ & 0.15052717690553594 & 0.33274431831757486\\
$x^2$ & $x+x^6$ & 0.01428216038702631 & 0.23635267804532455\\
$x^3$ & $x+x^6$ & -0.04777162399115301 & 0.17190364078803073\\
$x^4$ & $x+x^6$ & -0.08125607922853098 & 0.12268865661476232\\
$1$ & $x+x^2+x^3$ & 0.31281238144992430 & 0.42522475067652506\\
$x$ & $x+x^2+x^3$ & -0.03063713609272196 & 0.18754944674648485\\
$1$ & $2x+x^4$ & 0.22758105958079916 & 0.50285846929935148\\
$x$ & $2x+x^4$ & -0.08301824875147187 & 0.21632137918922486\\
$x^2$ & $2x+x^4$ & -0.12503111830114265 & 0.07854722450324651
\end{tabular}
\label{tab.Iinf}
\end{table}
\clearpage

Unfortunately some
fluky
maneouvre is still needed to decide on a transition
value
that divides the $x$-interval into the low region of good convergence of the sum $\sum p_{2k+1}$
and the high region of good convergence of the sum $\sum h_jX^j$, depending on
requirements
on
accuracy and number of terms of both series to be computed.
In practise one picks some $x$ larger than 1 where $q_{2k+1}$ converges,
computes the asymptotic series at the associate $X=1/x$, estimates
the error in the asymptotic part (which
adds some artistic flavor to the concept)
and
obtains a point spot estimate of $I_{p,\phi}(\infty)$ via (\ref{eq.Iinf}) subtracting both and
multiplying with $e^{i\phi(x)}$.
\begin{exa}
Continuing Example \ref{ex.asy}, one may compute a stable value
of $p(x)$ near $x\approx 1.8$ with 40 terms, illustrated in Figure \ref{fig.ada},
and compute an associate $p(X)$ again near $x\approx 1.8$ with the 
first three terms illustrated by the collocation of the $X^4$ and $X^5$
curves in Figure \ref{fig.asy}. The result is $I_{p,\phi}(\infty) \approx -0.07116+0.20606i$
with approximately two valid digits compared to Table \ref{tab.Iinf}. The imaginary
part matches the step height as $x\to \infty $ in Figure \ref{fig.asy}.
\end{exa}

\section{Summary}
The integral over the product of a polynomial $p(x)$ by a phasor $\exp[i\phi(x)]$ defined
by another polynomial $\phi(x)$ has been reduced to Incomplete Gamma Functions if $p$
and $\phi$ are powers of $x$, and to power series in $x$ and asymptotic series in $1/x$
for the generic case which differ by the general solution of an associated linear first-order
differential equation.

\appendix
\section{Reversion of basic polynomials} \label{sec.revSimpl}
Some elementary examples of series reversions (\ref{eq.xofy}) are gathered
in Table \ref{tab.revs}.
It shows the polynomial in $x$, the
initial few coefficients $\beta_j$ of its reversion,
and---if available---the label
of these coefficients
in the OEIS \cite{EIS}.
\begin{longtable}{lp{0.7\textwidth}p{0.1\textwidth}}
\caption{Basic Examples of Series reversions\label{tab.revs}}\\
$y=\ldots$ & $\beta_1$, $\beta_2$,$\beta_3$,\ldots & \cite{EIS} \\
\endfirsthead
\caption[]{Basic Examples of Series reversions (cont.)}\\
$y=\ldots$ & $\beta_1$, $\beta_2$,$\beta_3$,\ldots & \cite{EIS} \\
\endhead
\input{gentab.out}
\end{longtable}

To support economical evaluation of the reversions, the subsequent
list shows (homogeneous) hypergeometric recurrences
in annihilator notation.

These formats occur because the $\beta$ are determinantal forms
of the $\alpha$-coefficients of the polynomials \cite{Dominiciarxiv05,GesselEJC13,JanjicJIS15}.
From the point of view of Lagrange's inversion formula,
\begin{equation}
t! \beta_t = \frac{d^{t-1}}{dx^{t-1}}\left(\frac{x}{\sum_{k=1}^l \alpha_kx^k}\right)^t_{\mid x=0}
= \frac{d^{t-1}}{dx^{t-1}}\frac{1}{(\sum_{k=0}^{l-1} \alpha_{k+1}x^k)^t}_{\mid x=0}
\end{equation}
accompanying (\ref{eq.yofr}), evaluated with Fa\`a di Bruno's formula \cite{diBrunoASM6,CraikAMM112,JohnsonAMM109,TyrrellM9},
it portraits that the number of nonzero higher order derivatives of the inner function (the polynomial in
the denominator) is limited by the order $l$.
Expansion of the determinant along the (sparse) final row of
di-Bruno's determinant establishes the recurrence.

Two-term hypergeometric recurrences in this list lead to representations 
of $\beta$ as product-ratios of $\Gamma$-values, and to generating
functions for the $\beta$-sequence of the hypergeometric form.

The simplest group are the cases $y=x+kx^t$
with integer $k$ and $t$
where the coefficients $\beta$ are generalized Catalan numbers \cite{Copeland0512}:
\begin{itemize}
\item
The coefficients for series reversion of $y=x+kx^2$ are essentially
products of Catalan numbers with powers of $k$ and obey 2-term
recurrences.  Solving the quadratic equation for $x$ yields \cite{CallanJIS10}
\begin{multline}
x=[(1+4ky)^{1/2}-1]/(2k)
=\sum_{l\ge 1} \binom{1/2}{l}(4ky)^l/(2k)
\\
= -\sum_{l\ge 1} \binom{2l-2}{l-1}\frac{1}{l} k^{l-1}(-y)^l
.
\end{multline}
So $\beta_j = 
\binom{2j-2}{j-1}\frac{1}{j} (-k)^{j-1}
$
and $j\beta_j+2k(2j-3)\beta_{j-1}=0$.
\item
The non-vanishing coefficients of $y=x+x^3$ obey $(2k-1)(2k-2)\gamma_k
+3(3k-4)(3k-5)\gamma_{k-1} =0$,
where $\gamma_k=\beta_{2k-1}$.
\item
The non-vanishing coefficients of $y=x+x^4$ obey
$(3k-3)(3k-4)(3k-2)\gamma_k
+4(4k-5)(4k-7)(4k-6)\gamma_{k-1}
=0$, where $\gamma_k=\beta_{3k-2}$.
\item
The non-vanishing coefficients of $y=x+x^5$ obey
$(4k-4)(4k-5)(4k-3)(4k-6)\gamma_k
+5(5k-9)(5k-8)(5k-7)(5k-6)\gamma_{k-1}
=0$, where $\gamma_k=\beta_{4k-3}$.
\item
The non-vanishing coefficients of $y=x+2x^3$ obey
$(2k-1)(k-1)\gamma_k
+3(3k-4)(3k-5)\gamma_{k-1}
=0$, where $\gamma_k=\beta_{2k-1}$.
\item
The non-vanishing coefficients of $y=x+2x^4$ obey
$(3k-3)(3k-4)(3k-2)\gamma_k
+8(4k-5)(4k-7)(4k-6)\gamma_{k-1}
=0$, where $\gamma_k=\beta_{3k-2}$.
\end{itemize} 
Other recurrences are
\begin{itemize} 
\item
The coefficients of $y=x+x^2+x^4$ obey
$1147j(j-1)(j-2)\beta_j
+8(j-1)(j-2)(647j-738)\beta_{j-1}
+4(j-2)(224j^2+1504j-5157)\beta_{j-2}
+8(800j^3-5040j^2+8746j-2655)\beta_{j-3}
-192(4j-15)(2j-7)(4j-17)\beta_{j-4}
=0$.
\item
The coefficients of $y=x+x^3+x^4$ obey
$124j(j-1)(j-2)\beta_j
+(j-1)(j-2)(7j-88)\beta_{j-1}
+(j-2)(870j^2-3456j+3347)\beta_{j-2}
+(1243j^3-9870j^2+25869j-22490)\beta_{j-3}
+8(4j-15)(2j-7)(4j-17)\beta_{j-4}
=0$.
\item
The coefficients of $y=x+2x^2+x^3=x(1+x)^2$ obey
$2j(j-1)\beta_j
+3(3j-2)(3j-4)\beta_{j-1}
=0$, so $\beta_j = \prod_{i=2}^j (-3)(3i-2)(3i-4)/[2i(i-1)]
$.
\item
The coefficients of $y=x+x^2+2x^3$ obey
$7j(j-1)\beta_j
+16(j-1)(2j-3)\beta_{j-1}
+12(3j-1)(3j-7)\beta_{j-2}
=0$.

\end{itemize}

\section{Laurent Series} \label{app.asy}
The coefficients of the Laurent expansion (\ref{eq.asy})
are obtained by inserting the series and the polynomial $\phi$ (\ref{eq.phiofalpha}) into
(\ref{eq.Levin}):
\begin{equation}
X^{-m} = i q(X) \sum_{j=1}^l j\alpha_j X^{-j} -X^2 q'(X),
\end{equation}
where we assume that only the forms $p=x^m=1/X^m$, $m<l-1$, need to be studied because
the differential equation is linear and because higher powers $m$ have been decimated as proposed
in Remark \ref{rem.reduO}.
We compare coefficients of equal powers of $X$ in
\begin{equation}
X^{-m} = i \sum_{t=1}^\infty h_t X^ t \sum_{j=1}^l j\alpha_j X^{-j} -X^2 \sum_{t=1}^\infty th_tX^{t-1}.
\end{equation}
Elementary resummation of the double sum on the right hand side suggests
\begin{equation}
X^{-m} = i \sum_{s=-l}^0 \sum_{j=-s+1}^l j\alpha_j h_{s+j} X^s
+ i \sum_{s=1}^\infty \sum_{j=1}^l j\alpha_j h_{s+j} X^s -\sum_{s=1}^\infty (s-1)h_{s-1}X^s.
\label{eq.solvasy}
\end{equation}
For $s\le 0$, only the first double sum on the right hand side contributes, so the
$h_t$ are obtained recursively for increasing $t=1,2,\ldots ,l$ via
\begin{equation}
i l\alpha_l h_t =
\delta_{t-l,-m}
- i \sum_{j=-t+1}^{-1} (j+l)\alpha_{j+l} h_{t+j},\quad 1\le t \le l.
\end{equation}
For $s> 0$, the left hand side of (\ref{eq.solvasy}) vanishes and only the second and third sum on the right hand side contribute,
so the
$h_t$ are obtained recursively for $t=l+1, l+2,\ldots$ via
\begin{equation}
il\alpha_l h_t = (t-l-1)h_{t-l-1}- i \sum_{j=-l+1}^{-1} (j+l)\alpha_{j+l} h_{t+j} ,\quad l<t.
\end{equation}
with the convention that $h_t=0$ whenever $t\le 0$.

\bibliographystyle{amsplain}
\bibliography{all}

\end{document}